\documentclass[11pt,a4paper,leqno]{article}
\usepackage{amsmath,amsthm,amsfonts,amssymb,amscd}
\usepackage{graphics}
\usepackage[latin1]{inputenc}

\def \é{\'e}
\def \è{\`e}
\def \à{\`a}
\def \ê{\^e}
\def \ù{\`u}
\def \î{\^{\i}}
\def \â{\^a}
\def \û{\^u}
\def \ô{\^o}
\def \ï{\"{\i}}
\def \ö{\"{o}}
\def \ç{\c{c}}
\def\unt{{\underline\theta}}

\def \te {\theta}

\newcommand{\la}{{\lambda}}
\newcommand{\ov}{\overline}
\newcommand{\re}{{\mathbb{R}}}
\newcommand{\al}{{\alpha}}

\newcommand{\ve}{{\varepsilon}}

\newcommand{\und}{\underline}
\newcommand{\teta}{\theta}

\def\cb{{\vskip-4mm{\hfill\hbox{$\vrule\vcenter to 2mm{\hrule\vfil\hbox{\kern
2mm}\vfil\hrule}\vrule$}\quad}} \medskip}

\begin{document}

\centerline{\Large\bf Geometric properties of images of cartesian products}

\vskip .1in

\centerline{\Large\bf of regular Cantor sets by differentiable real maps}

\vskip .2in

\centerline{\sc Carlos Gustavo Moreira}

\centerline{\sl Instituto de Matem\'atica Pura e Aplicada}

\centerline{\sl Estrada Dona Castorina 110}

\centerline{\sl 22460-320 Rio de Janeiro, RJ}

\centerline{\sl Brasil}
\vskip .2 in
{\small \hskip 2.4in \it Dedicated to Jean-Christophe Yoccoz and Welington de Melo}

\vskip .2in
{\bf Abstract:} We prove dimension formulas for arihmetic sums of regular Cantor sets, and, more generally, for images of cartesian products of regular Cantor sets by differentiable real maps.
\vskip .2in
\noindent
{\bf I. Introduction:} 

We prove dimension formulas for arihmetic sums of regular Cantor sets, and, more generally, for images of cartesian products of regular Cantor sets by differentiable real maps.

We say that a $C^2-$ regular Cantor set $K$ is {\it non essentially affine} if the expanding dynamics which defines it cannot be $C^2$ conjugated to a map whose second derivative is identically zero on the Cantor set. We have the following result:

\noindent{\bf Theorem 1}: Let $K$\,, $K'$ be regular Cantor sets with Hausdorff dimensions $d$ and $d'$, respectively, such that $K$ is of class $C^2$ and non essentially affine.
Then, given a $C^1$ map $f$ from a neighbourhood of $K \times K'$ to $\re$ such that, in some point of $K \times K'$ its gradient is not parallel to any of the two coordinate axis, we have $HD(f(K \times K'))=\min\{1, d+d'\}$. where $HD$ denotes Hausdorff dimension.
\vglue .1in
\noindent{\bf Remark}: Since $K$ and $K'$ are regular Cantor sets, we have \hfill\break $\ov{\dim}(K)=HD(K)=d$, $\ov{\dim}(K')=HD(K')=d'$ and $\ov{\dim}(K \times K')=HD(K \times K')=d+d'$,  where $\ov{\dim}$ denotes upper box dimension, and so, since $f$ is Lipschitz, \hfill\break $\ov{\dim}(f(K \times K'))\le \min\{1, \ov{\dim}(K \times K')\}=\min\{1, d+d'\}=HD(f(K \times K'))\le \ov{\dim}(f(K \times K'))$, therefore we also have $\ov{\dim}(f(K \times K'))=\min\{1, d+d'\}$.

The first version of these formula were stated in [Mor98], without a proof, as a tool to obtain results on geometric properties of the classical Markov and Lagrange spectra. Shmerkin wrote in 
[Shm] a proof of Theorem 1 above assuming that both $K$ and $K'$ are $C^2$ and a non-resonance hypothesis: that there are periodic orbits of the dynamics which define $K$ and $K'$ such that the ratio of the logarithms of the absolute values of their eigenvalues is irrational. 

An immediate corollary is the following:

\noindent{\bf Proposition}: Let $K$\,, $K'$ be regular Cantor sets
with Hausdorff dimensions $d$ and $d'$, respectively, such that $K$ is of class $C^2$ and non essentially affine. Then $HD(K+sK') =\ov{\dim}(K+sK')=
\min\{1, HD(K) + HD(K')\}$, for every $s \neq 0$.

This formula is an essential tool for the proof of the main result of the paper [M] on geometric properties of the classical Markov and Lagrange spectra.

We also prove a version of these results in the case when $K$ and $K'$ are not necessarily $C^2$, assuming that they satisfy the hypothesis of the Scale Recurrence Lemma of [MY2], and that they have periodic orbits whose ratio of logarithms of (the norms of the) eigenvalues is irrational. This result is an essential tool in the work [CMMR] in collaboration with Cerqueira and Matheus on geometric properties of dynamical Markov and Lagrange spectra associated to conservative diffeomorphisms and in a forthcoming work with Cerqueira, Matheus and Roma\~na on geometric properties of dynamical Markov and Lagrange spectra associated to geodesic flows in negative curvature.

I would like to thank Aline Cerqueira, Carlos Matheus and Pablo Shmerkin for helpful comments and suggestions which substantially improved this work.

\vskip .1in

\noindent

\vskip 0,3cm
{\bf 1 - Basic definitions}.

{\bf 1.1 - Regular Cantor sets.}

\vskip 0,3cm

We will follow the notations and definitions of [MY] and [MY2]. 

Let ${\bf A}$ be a finite alphabet and $\Sigma \subset {\bf A}^{\mathbb Z}$ a
subshift of finite type defined by a transition set ${\cal B} \subset {\bf A}^2$. We
will always assume that $\Sigma$ is mixing and uses all letters of ${\bf A}$.

A {\bf word} of $\Sigma$ is a finite sequence of elements of A verifying the transition rules given by ${\cal B}$.

Let $r \in (1, + \infty]$ ; we say that a $C^r$ map $g$ is
{\bf expanding of type $\Sigma$} if:

- its domain is a finite union $\displaystyle \bigsqcup_{\cal B} \;
I(a_0, a_1)$, where, for every $(a_0 , a_1) \in {\cal B} \; , \; I(a_0, a_1)$ is a
compact sub-interval of $I (a_0):=[0,1] \times \{a_0\}$ ;

- for each $(a_0 , a_1) \in {\cal B}$, the restriction of $g$ to $I(a_0 , a_1)$ is a
$C^r$-diffeomorphism onto $I(a_1)$, satisfying $\vert Dg(t)\vert >1$ for every
$t \in I(a_0 , a_1)$.

We denote by $\Omega^r_\Sigma$ the space  $C^r$ expanding maps of type
$\Sigma$, endowed with the $C^r$-topology. Such a map defines the {\bf
regular Cantor set}

$$K = \bigcap_{n \geq 0} \; g^{-n} (\; \bigsqcup_{\cal B} \; I(a_0, a_1)).$$

The dynamics of $g$ over $K$ is canonically conjugated to the unilateral sub-shift
$\sigma:\Sigma^+\to \Sigma^+$, $\Sigma^+ \subset {\bf A}^{\mathbb N}$ defined by ${\cal B}$: there exists a unique homeomorphism $h\colon \Sigma^+ \to K$ such that
\begin{align*}
&h(\underline{a}) \in I(a_0), \text{ for } \underline{a} = (a_0,a_1,\dots) \in \Sigma^+\\
&h\circ\sigma = g \circ h.
\end{align*}

\vskip 0,5cm
{\bf 1.2 - Limit Geometries.}

For $(a_0, a_1) \in {\cal B} , f_{(a_0, a_1)}$ denotes the inverse map of $g|_{I(a_0 ,
a_1)}$. For a word $\underline a = (a_0, \dots , a_n)$ of $\Sigma$, we put
$$f_{\underline a} = f_{a_0 a_1} \circ \dots \circ f_{a_{n-1}a_n} \; ;$$ it is a
diffeomorphism of $I(a_n)$ over a sub-interval $I(\underline a)$ of $I(a_0).$

We denote $K(\theta_0):=K \cap I(\theta_0)$ and, more generally, $K(\und a):=K \cap I(\und a)$.

We denote by $\Sigma^-$ the unilateral sub-shift defined by $\cal B$ indexed by the integers which are negative or zero.

We equip $\Sigma^-$ with the following ultrametric distance: for ${\und\teta} \ne \tilde{\und\teta} \in \Sigma^-$, we set 
$$
d(\und{\teta},\und{\widetilde\teta}) =
\begin{cases}
1 \quad\text{if}\quad \teta_0 \ne \widetilde\teta_0\\
|I(\und{\teta} \wedge \und{\widetilde\teta})|\quad\text{otherwise}
\end{cases}
$$
where $\und{\teta} \wedge \und{\widetilde\teta} = (\teta_{-n},\dots,\teta_0)$ with
$\widetilde\teta_{-j} = \teta_{-j}$ for $0 \le j \le n$ and $\widetilde\teta_{-n-1}
\ne \teta_{-n-1}$\,.

For ${\underline \theta} = (\theta_m)_{m \leq 0}$ in
$\Sigma^- \; , \; n
< 0$, we put
$$k^{\underline \theta}_n = B \circ f_{\underline \theta^n} \; , $$ where
${\underline \theta}^n
= (\theta_{-n}, \dots , \theta_0)$ and $B$ is the only affine map from
$I({\underline
\theta}^n)$ onto $I({\theta}_0)$ such that $k_n^{\underline \theta}$
is orientation preserving.

The sequence $k^{\underline \theta}_n$ converges to a diffeomorphism
$k^{\underline \theta} \in Diff_+^r(I(\theta_0))\text{ (cf. [Su])}$ ; the convergence 
takes place in the $C^{r'}$-topology for every $r' <r$, and even in the $C^r$-topology if $r$
is an integer or $+\infty$. The convergence is uniform in
${\underline \theta} \in \Sigma^-$ and in a neighbourhood of $g$ in
$\Omega^r_\Sigma$.

The {\bf limit geometry} of $K$ associated to ${\underline \theta}$
is the Cantor set
$$K^{\underline \theta} = k^{\underline \theta} (K \cap I(\theta_0)).$$

The diffeomorphisms $k^{\underline \theta} , {\underline \theta} \in
\Sigma^-$, realise a fibered linearisation of the dynamics since the applications
$k^{\sigma^{-1}(\underline \theta)} \circ f_{\theta_{-1}\theta_0} \circ
(k^{\underline\theta})^{-1}$ are affine.

We denote $I^{\und \theta}(\und a):=k^{\underline \theta}(I(\und a))$.

The different $K^{\und{\te}}$ are related in the following way: let
$F^{\und{\te}}$ be the affine map from $I(\te_0)$ onto $I^{\sigma^{-
1}(\und{\te})}(\te_{-1},\te_0)$ with the same orientation than 
$f_{\te_{-1}\te_0}$\,. Then, on $I(\te_0)$:
$$
F^{\und{\te}} \circ k^{\und{\te}} = k^{\sigma^{-1}(\und{\te})} \circ f_{\te_{-1}\te_0}\,.
$$
It follows by induction that, for any positive integer $n$, if we denote by $F_n^{\und{\te}}$ the affine map from $I(\te_0)$ onto $I^{\sigma^{-n}(\und{\te})}({\underline \theta}^n)$ with the same orientation than $f_{{\underline \theta}^n}$, then, on $I(\te_0)$:
$$
F_n^{\und{\te}} \circ k^{\und{\te}} = k^{\sigma^{-n}(\und{\te})} \circ f_{{\underline \theta}^n}\,.
$$
Taking inverses, we also get that $k^{\und{\te}} \circ (g^n|_{I({\underline \theta}^n)})\circ (k^{\sigma^{-n}(\und{\te})})^{-1}=(F_n^{\und{\te}})^{-1}$ is an affine map from $I^{\sigma^{-n}(\und{\te})}({\underline \theta}^n)$ onto $I(\te_0)$.
\vskip 0,5cm

{\bf 1.3 - Renormalization Operators.}

For $a \in {\bf A}$, we denote by ${\cal P}^r(a)$ the space of
$C^r$-embeddings of $I(a)$
into ${\re}$, endowed with the $C^r$-topology ; we put ${\cal P}^r =
\displaystyle\bigsqcup_{\bf A}
{\cal P}^r(a)$. We call by ${\cal A}$ the space of pairs $({\underline \theta},
A)$, where
${\underline \theta} \in \Sigma^-$ and $A$ is an {\bf affine} embedding of
$I(\theta_0)$ into ${\re}$ ; we have a natural map: ${\cal A}
\rightarrow {\cal P}^r$ which associates $A \circ k^{\underline \theta}$ to the pair $({\underline
\theta}, A)$.

Let $\underline a = (a_0, \dots , a_n)$ a word of $\Sigma$ ; the {\bf renormalization operator} $T_{\underline a} : {\cal P}^r(a_0) \rightarrow {\cal P}^r(a_n)$ is defined by 
$$T_{\underline a} (h)= h \circ f_{\underline a} \; .$$

It has a lift to ${\cal A}$ (still called by $T_{\underline a}$) ;
while
${\underline a}$ varies between the words of $\Sigma$ of size $n+1$, the
$T_{\underline a}$ are the inverse branches of $S^n$, the map $S
: {\cal A}
\rightarrow {\cal A}$ being defined by

$$S({\underline \theta}, A)= ( \sigma^{-1} {\underline \theta} , A \circ
k^{\underline
\theta} \circ g|_{I(\theta_{-1} , \theta_0)} \circ (k^{\sigma^{-1}({\underline
\theta})})^{-1}) \; . $$

{\bf 1.4 - Relative Configurations.}

As in 1.3, we can introduce the spaces $\cal A$, $\cal A'$ associated to the limit geometries of $g$,\,\,$g'$ respectively. We denote by $\cal C$ the quotient of ${\cal A}\times{\cal A}'$ by the diagonal action on the left of the affine group. An element of $\cal C$, represented by $(\underline\theta,A) \in \cal A$,\,\, $(\underline{\theta}', A') \in \cal A'$, is called a relative configuration of the limit geometries determined by $\underline{\theta}$, $\underline{\theta}'$.

The fibers of the quotient map: ${\cal C} \to {\cal S}$ are 
one-dimensional and have a canonical affine structure; moreover this bundle
map is trivializable; we choose now an explicit trivialization ${\cal C} \cong {\cal S}
\times {\mathbb R}$ in order to have a coordinate in each fiber.

For each $a \in {\mathbb A}$ (resp. $a' \in {\mathbb A}'$), we choose $\omega(a) \in \Sigma^+$
starting with $a$ (resp. $\omega(a') \in \Sigma^{\prime +}$ starting with $a'$).
Given a configuration represented by $(\underline{\theta},A)$, $(\underline{\theta}',A')$, we
normalize (by the action of the affine group) in order to obtain that $A'\colon
{\mathbb R}\times\{\theta_0\} \to {\mathbb R}$ becomes a translation sending
$k^{\prime\underline{\theta}'}(h'(\omega(\theta_0')))$ to 0.

Then our coordinate $t((\underline{\theta},A), (\underline{\theta}',A'))$ on the fiber will be by
definition
$$
t = A(k^{\underline{\theta}}(h(\omega(\theta_0)))).
$$

With the chosen normalization, $A(I(\theta_0))$ has size 1 while $A'(I(\theta_0'))$
has size $|s|$.

\vskip 0,5cm
{\bf 2 - The scale recurrence property}.

{\bf 2.1} - Let $g$, $g'$ be expansive maps as above. If
$((\und{\te},A), (\und{\te}',A'))$ is a relative configuration, the affine map $A^{-1}\circ A'\colon \re\times\{\te_0'\} \to \re\times\{\te_0\}$ is well defined; the coefficient of $x$ is a number in $\re^*$ which is called the {\it relative scale\/} of the given relative configuration. We obtain thus the space
${\cal S} = \Sigma^- \times \Sigma^{\prime -} \times \re^*$ of relative scales;
this space is
considered as a quotient of ${\cal C}$ in the following way: for
$(\underline \theta , A) \in {\cal A}\; , \; (\underline \theta' , A') \in {\cal
A'}$, we put
$s = \frac{DA'}{DA}$ ; the map
$$((\underline \theta , A), (\underline \theta' , A')) \longmapsto  (\underline
\theta \, ,
\, \underline \theta' ,s)$$

\noindent defines by passing to the quotient a projection from ${\cal C}$ to
${\cal S}$.

Given words $\underline{a} =(a_0,a_1,\dots,a_n)$ and $\underline{a}' =(a_0',a_1',\dots,a_n')$, we define renormalization operators $T_{\und{a}}$, $T'_{\und{a}'}$ acting on ${\cal S}$ in the following way:

We have, for $s \in \re^*$
$$
T_{\und{a}}(\und{\te}, \und{\te}',s) = (\und{\te}a_1\dots a_n, \und{\te}', \ve
s |I^{\und{\te}}(\und{a})|^{-1})=(\und{\te}\sigma(\und{a}), \und{\te}', \ve
s |I^{\und{\te}}(\und{a})|^{-1})
$$
where $\sigma(\und{a})=\sigma(a_0,a_1,\dots,a_n)=(a_1,\dots,a_n)$ and $\ve = +1$ (resp. $-1$) if $f_{\und{a}}$ is orientation preserving (resp.
reversing).

Similarly, if $\und{a}'$ is a word in $\Sigma$ starting with $\te_0'$\,, we
have
$$
T'_{\und{a}'}(\und{\te},\und{\te}',s) = (\und{\te}, \und{\te}' a_1'\dots a_n',
\ve's |I^{\und{\te}'}(\und{a}')|)=(\und{\te}, \und{\te}' \sigma(\und{a}'),
\ve's |I^{\und{\te}'}(\und{a}')|).
$$

For $R > 1$, we put $J_R = [-R , -R^{-1}] \cup [R^{-1}, R]$; we define ${\cal S}_R = \Sigma^- \times \Sigma^{'-} \times J_R$.

\vskip 0,5cm
{\bf 2.2} - The formulation of the scale recurrence property involves, between several other constants, a parameter $\alpha \in (0,1]$ which will play an important rôle in what follows. In the case of [MY] (which corresponds to the case when both Cantor sets are $C^2$), we had $\alpha = 1$. In the case of Cantor sets arising as stable and unstable Cantor sets associated to a horseshoe of a diffeomorphism, we may have to take $\alpha<1$, related to the contraction and expansion rates of (iterates of) the diffeomorphism near the horseshoe, which also affects the differentiability classes of the corresponding stable and unstable Cantor sets. On the other hand we state the scale recurrence property and the scale recurrence lemma in terms of regular Cantor sets (not a priori associated to a horseshoe).
 
\vskip 0,5cm
{\bf 2.3} - We choose once for all a constant $c_0 > 1$. For every $0 <
\rho < 1$, we then denote by $\Sigma(\rho)$ the set of words $\underline a$ of
$\Sigma$ such that $c^{-1}_0 \rho \leq \vert I( \underline a ) \vert \leq c_0 \rho $. We define analogously $\Sigma'(\rho)$.

The scale recurrence property can now be stated as follows:

The constants $c_0 , R$ being choosed large enough, there are $c_1 ,
c_2, c_3 > 0 ,
\rho_0 \in (0,1)$ such that, for every $0 < \rho < \rho_0$, and for every collection of sets $(E(\underline a , \underline a'))_{(\underline a , \underline a')
\in
\Sigma(\rho) \times \Sigma'(\rho)}$ verifying
$$E(\underline a , \underline a') \subset J_R \; ,$$
$$Leb(J_R - E(\underline a , \underline a')) < c_1 \; ,$$ we can find a
collection of compact sets $(E^*(\underline a , \underline a'))_{\underline a ,
\underline a' \in \Sigma(\rho) \times \Sigma'(\rho)}$ satisfying

(i) $E^*(\underline a , \underline a')$ is contained in the $c_2
\rho^\alpha$-neighbourhood of $E(\underline a , \underline a')$ in $J_R$ ;

(ii) For at least half of the pairs $(\underline a , \underline a')$, the Lebesgue measure of 
$E^*(\underline a , \underline a')$ is larger than half of that of $J_R$ ;
 
(iii) for each $(\underline a , \underline a') \in \Sigma(\rho) \times
\Sigma'(\rho)$, and each $s \in E^*(\underline a , \underline a')$, there are at least $c_3
\rho^{-(d+d')}$ pairs $(\underline b , \underline b') \in \Sigma(\rho)
\times \Sigma'(\rho)$ so that we have for each $\underline \theta \in
\Sigma^- ,  \underline \theta' \in \Sigma^{'-}$ ending respectively by
$\underline a , \underline a'$,
$$T_{\underline b} \; T'_{\underline b'} (\underline \theta , \underline
\theta', s)= (\tilde
{\underline \theta} , \tilde {\underline \theta'}, \tilde s)$$ with $[\tilde s - \rho^\alpha , \; \tilde s
+
\rho^\alpha] \subset E^*(\underline b \; ,
\; \underline b').$

{\bf Remark} - The pairs $(\underline b ,
\underline b')$ appearing in (iii) should of course have as first letters the last letters of $(\underline a , \underline a')$.

{\bf 2.4} - We proved in [MY2] that the scale recurrence property is a consequence of the following property $(H \alpha)$, which deals only with the first Cantor set $K$.

{\bf Property $(H \alpha)$} : There exist $\eta >0 \; , \; \rho_1 \in (0,1)$
such that, for every $0<\rho < \rho_1 \; , \; 1 \leq \xi \leq \rho^{-\alpha}$ and every subset
$X$ of $\Sigma(\rho)$ satisfying
$$\#(\Sigma(\rho) - X) < \eta \rho^{-d} \; ,$$ we can find $\underline a^0 ,
\underline
a^1 \in X$ such that we have, for every $\Phi
\in {\re}$ and every $\underline \theta^0 \in \Sigma^- \; , \; \underline
\theta^1 \in
\Sigma^{'-}$ ending by $\underline a^0 , \underline a^1$,
$$
\# \bigg\{\underline b \in \Sigma(\rho), \bigg\vert \; \sin \; \bigg[ \frac{1}{2} \;
\xi
\;
\log  \frac{\vert I^{\underline \theta^0} (\underline b) \vert}{\vert
I^{\underline
\theta^1} (\underline b) \vert} + \Phi \bigg] \; \bigg\vert \;\geq \eta \bigg\}
\geq
\eta  \rho^{-d} \; . $$

When $\alpha=1$, and so the $k^{\underline \theta}$ are $C^2$, $(H1)$ follows from the much simpler property

\noindent
$(H'1)$ \qquad There are $\underline \theta^0 , \underline \theta^1 \in
\Sigma^-$, ending by the same letter $\theta_0$, and $x_0 \in I(\theta_0)$ such that we have
$$D \; \log   D \; [k^{\underline \theta^1} \circ (k^{\underline \theta^0})^{-1}
] \; (x_0)
\not= 0 \; .$$

When $K$ and $K'$ come from horseshoes $\Lambda$ and
$\Lambda'$ associated to a diffeomorphism $F$, i.e, are, respectively, a stable Cantor set of $\Lambda$ and an unstable Cantor set of $\Lambda'$, we choose $\alpha$ in the following way.

When the diffeomorphism $F$ contracts area in a neighbourhood of the horseshoe $\Lambda$ and expands area in a neighbourhood of the horseshoe $\Lambda'$, we choose $\alpha=1$.

In the conservative case (when $F$ preserves a smooth measure) we choose arbitrarily $\alpha \in (0,1)$.

In the remaining cases, we choose $\alpha \in (0,1)$ in such a way that the following property is satisfyied: for an appropriate riemannian metric on $M$, we have
$$\Vert T_z F_0 |_{E^s} \Vert < \Vert T_z F_0|_{E^u} \Vert^{-\alpha}\; ,$$
$$\Vert T_{z'} F_0^{-1} |_{E^u} \Vert < \Vert T_{z'} F_0^{-1}|_{E^s}
\Vert^{-\alpha}\; ,$$
for every $z \in \Lambda \; , \; z' \in \Lambda' \; . $

In all cases, this assures that there exists $\alpha' > \alpha$ such that the foliations
$W^s(\Lambda)$ and $W^u(\Lambda')$ are transversely of class
$C^{1+\alpha'}$. The linearizing maps $k^{\underline \theta} ,
k^{'\underline
\theta'}$ will be of class $C^{1+\alpha'}$ as well. In this case the correspondences $\underline \theta \to k^{\underline \theta}$, $\underline{\theta}'\to k'^{\underline{\theta}'}$ are $\alpha$-Hölder.

{\bf Remark} - The hypothesis $(H \al)$ are less restrictive when $\al$ is smaller. So, when $K$ is $C^2$ (even when $K'$ is not $C^2$) (in particular, when $K$ and $K'$ come from horseshoes as before, $F$ contracts area in a neighbourhood of the horseshoe $\Lambda$ but does not expand area in a neighbourhood of the horseshoe $\Lambda'$) it is enough to check the property $(H'1)$, which makes the situation in this case simpler than in the general case.

It is proved in [MY2] that, for Cantor sets coming from horseshoes, the properties $(H \alpha)$ (or $(H'1)$ when $\alpha =1$) are generic (open and dense) conditions in both conservative and dissipative cases, and so is the scale recurrence property. Since the properties $(H \alpha)$ follows from the property $(H'1)$, this is also true for $C^k$ Cantor sets, $k\ge 2$.

\vskip 0,5cm
{\bf 3 - Dimension estimates - a discrete Marstrand-like result}.

Let $K,K'$ be regular Cantor sets with associated finite alphabets ${\bf A}$, ${\bf A'}$. Let $\tilde K$ and $\tilde K'$ be differentiable embeddings of $K$ and $K'$, respectively, i.e., $\tilde K=h(K\cap I(a))$, $\tilde K'=\tilde h(K'\cap I'(a')))$, for some $C^1$-embeddings $h\in {\cal P}^1(a)$, $\tilde h\in {\cal P}^1(a')$, with $a\in {\bf A}$, $a'\in {\bf A'}$. Given a word $\und c$ of $\Sigma$ beginning by $a$, we define the interval $\tilde I(\und c):=h(I(\und c))$. Analogously, given a word $\und c'$ of $\Sigma'$ beginning by $a'$, we define the interval $\tilde I(\und c'):=\tilde h(I(\und c'))$. Given $\delta \in (0,1)$, we denote by $\tilde \Sigma(\delta)$ the set of words $\und c\in \Sigma(\delta)$ which begin by $a$, and by $\tilde \Sigma'(\delta)$ the set of words $\und c'\in \Sigma'(\delta)$ which begin by $a'$. We say that $Q$ is a {\it $\delta$-rectangle} of $\tilde K \times \tilde K'$ if there is $(\und c, \und c') \in \tilde \Sigma(\delta) \times \tilde \Sigma'(\delta)$ such that $Q=\tilde I(\und c) \times \tilde I'(\und c')$. Let ${\cal R}(\delta)=\tilde \Sigma(\delta) \times \tilde \Sigma'(\delta)$ be the set of $\delta$-rectangles of $\tilde K \times \tilde K'$.

\noindent
{\bf Proposition 1:} {\it Let $K,K'$ be regular Cantor sets with $d+d'<1$, where $d$ and $d'$ are the Hausdorff dimensions of $K$ and $K'$, respectively, and let $\tilde K$ and $\tilde K'$ be differentiable embeddings of $K$ and $K'$, as above. For $s\in \mathbb R$, let $\pi_s:{\mathbb R}^2\to \mathbb R$ be given by $\pi_s(x,y)=x-sy$. Let $R>0$ be a large constant. For a given $\delta>0$, let $N_{\delta}(s)$ be given by $N_{\delta}(s)=\#\{(Q,\tilde Q) \in {\cal R}(\delta) \times {\cal R}(\delta) \mid \pi_s(Q)\cap\pi_s(\tilde Q)\ne \emptyset\}$.

Then $\int_{-R}^RN_{\rho}(\la)d\la=O(\rho^{-d-d'})$}, $\forall \rho \in (0,1)$.

\vskip .1in

\noindent
{\bf Proof:} Let $m=\lfloor\log_2\rho^{-1}\rfloor$. For $0\le r \le m$, let us estimate the measure of the set $\{\la\in[-R,R]\mid \pi_\la(Q)\cap\pi_\la(\tilde Q)\ne \emptyset\}$, where $Q,\tilde Q$ are $\rho$-rectangles of $\tilde K \times \tilde K'$ with $2^{-r}\le d(Q,\tilde Q)< 2^{-r+1}$. This set is an interval. Let $p_1=(x_1,y_1)$ and $p_2=(x_2,y_2)$ be the centers of $Q$ and $\tilde Q$, respectively. If $\pi_\la(Q)\cap\pi_\la(\tilde Q)\ne\emptyset$, we should have $d(\pi_\la(p_1),\pi_\la(p_2))=O(\rho)$, i.e., $|(x_1-\la y_1)-(x_2-\la y_2)|=O(\rho)$, so
$$
\left\vert \la - \frac{x_1-x_2}{y_1-y_2} \right\vert= O\left( \frac{\rho}{|y_1-y_2|} \right).
$$
If $|y_1-y_2|<2^{-r-2}/R$, $|x_1-x_2|> 2^{-r-1}$, so
$$
|(x_1-\la y_1)-(x_2-\la y_2)|=|(x_1-x_2)-\la(y_1-y_2)|>2^{-r-1}-2^{-r-2}=2^{-r-2}>\rho
$$
for $\la\in[-R,R]$ and $r\le m-2$. If $|y_1-y_2|\ge 2^{-r-2}/R$ then $O\left( \frac{\rho}{|y_1-y_2|}\right)=O(2^r\rho)$, so
$$
m(\{\la\in[-R,R] \mid \pi_\la(Q)\cap\pi_\la(\tilde Q)\ne\emptyset\})=O(2^r\rho).
$$

Let us now estimate $\int_{-R}^R N_{\rho}(\la) d\la$, where $N_{\rho}(\la)$ is given by $N_{\rho}(\la)=\#\{(Q,\tilde Q)$ $\rho$-rectangles of $\tilde K \times \tilde K' \mid \pi_\la(Q)\cap\pi_\la(\tilde Q)\ne \emptyset\}$. If we fix a $\rho$-rectangle $Q$ of $\tilde K \times \tilde K'$, the number of $\rho$-rectangles $\tilde Q$ of $\tilde K \times \tilde K'$ with $2^{-r}\le d(Q,\tilde Q)<2^{-r+1}$ is $O((2^{-r})^{d+d'}/\rho^{d+d'})$, and the number of $\rho$-rectangles $Q$ of $\tilde K \times \tilde K'$ is $O(\rho^{-d-d'})$. So we have:
\begin{align*}
\int_{-R}^R N_{\rho}(\la) d\la &= \sum_{Q\in{\cal R}(\rho)} \sum_{r=0}^m \sum_{\substack{\tilde Q\in{\cal R}(\rho) \\ 2^{-r}\le d(Q,\tilde Q)<2^{-r+1}}} m(\{\la\in[-R,R] \mid \pi_\la(Q) \cap \pi_\la(\tilde Q)\ne \emptyset\})= \\
&= O(\rho^{-d-d'} \sum_{r=0}^m(2^{-r}\rho^{-1})^{d+d'} \cdot 2^r\rho)=O(\rho^{1-2d-2d'} \sum_{r=0}^m 2^{r(1-d-d')})= \\
&= O(\rho^{1-2d-2d'}\cdot 2^{m(1-d-d')}) = O(\rho^{1-2d-2d'}(\rho^{-1})^{1-d-d'})=O(\rho^{-d-d'})
\end{align*}
 \qed
\noindent

Given a point $(\underline{\theta}, \underline{\theta}',s)$ in the space
$\cal S$ of relative scales, and points $x \in K(\theta_0)$, $x' \in K'(\theta_0')$,
there is a unique relative configuration above $(\underline{\theta}, \underline{\theta}',s)$ such
that
$$
A(k^{\underline{\theta}}(x)) = A'(k^{\prime\underline{\theta}'}(x')),
$$
(where $(\underline{\theta},A), (\underline{\theta}',A')$ represents this configuration).
Parametrizing the fiber of $\cal C$ over $(\underline{\theta}, \underline{\theta}',s)$ as in \S1.4,
the coordinate of this configuration is given by
\begin{align*}
t &:= \pi_{\underline{\theta},\underline{\theta}',s}(x,x')\\
&= k^{\underline{\theta}}(x) - k^{\underline{\theta}}(h(\omega(\theta_0))) -
s(k^{\prime\underline{\theta}'}(x') - k^{\prime\underline{\theta}'}(h(\omega(\theta_0')))).
\end{align*}

Notice that if $h(x)=k^{\underline{\theta}}(x) - k^{\underline{\theta}}(h(\omega(\theta_0)))$ and $\tilde h(x')=k^{\prime\underline{\theta}'}(x') - k^{\prime\underline{\theta}'}(h(\omega(\theta_0')))$, then $\tilde K=h(K\cap I(\theta_0))$, $\tilde K'=\tilde h(K'\cap I'(\theta_0')))$ are affine embeddings of $K^{\underline\theta}$ and $K^{'\underline\theta'}$, respectively, and $\pi_{\underline \theta , \underline \theta', s} (x,x')=\pi_s(h(x), \tilde h(x'))$ is the projection by $\pi_s$ of the point $(h(x), \tilde h(x'))\in \tilde K \times \tilde K'$.

For a given $\delta>0$, let $N_{\delta}(\underline \theta , \underline \theta', s)$ be given by $N_{\delta}(\underline \theta , \underline \theta', s)=\#\{((\und a, \und a'),(\und{\tilde a},\und{\tilde a'})) \in (\Sigma(\delta) \times \Sigma'(\delta))^2 \mid \pi_{\underline \theta , \underline \theta', s}(I_{\und a} \times I_{\und a'}')\cap\pi_{\underline \theta , \underline \theta', s}(I_{\tilde{\und a}} \times I_{\tilde{\und a'}}')\ne \emptyset\}$.

{\bf Remark 1:} 
 It follows from the proof that the estimate in the $O(\rho^{-d-d'})$ notation above depend continuously on $K,K',R$ (and on $1-d-d'$). In particular, there is a constant $c_4>0$ such that we have $\int_{-R}^RN_{\delta}(\underline \theta , \underline \theta', \la)d\la \le c_4 \cdot \delta^{-d-d'}$ for every $\delta\in (0,1]$ and every pair $(\unt,\unt')$ of limit geometries of $K$ and $K'$.
\vskip .1in
{\bf Remark 2:} 
Notice that the constant implicit in the $O$ notation gets worst when $d+d'$ approaches 1, which is consistent with Besicovitch's theorem on projections of irregular 1-sets.

Proposition 1 is useful to estimate measures and dimensions of projections because of the folowing result.

{\bf Proposition 2:} Let $\delta \in (0,1]$. Suppose that $d+d'<1$ and $N_{\delta}(\la)\le \tilde c \cdot \delta^{-d-d'}$, where $\tilde c>0$ is a constant. Let ${\cal F} \subset {\cal R}(\delta)$ be a family of $\delta$-rectangles of $\tilde K \times \tilde K'$ such that $\#{\cal F}\ge b \cdot \delta^{-d-d'}$, where $b>0$ is a constant. Let $X=(\tilde K \times \tilde K')\cap(\cup_{Q \in {\cal F}}Q)$. Then we have $Leb(\pi_{\la}(X))\ge c_0^{-1} \cdot {\tilde c}^{-1} \cdot b^2\cdot \delta^{1-d-d'}/4.$ 

{\bf Proof:} Let $k=\lfloor (2R+2)c_0\delta^{-1} \rfloor$, and let $x_0=-R-1, x_1, x_2,\dots, x_k, x_{k+1}=R+1\in [-(R+1),R+1]$ be an arithmetic progression. Since $\pi_\la(Q)$ has size at least $c_0^{-1}\delta>(2R+2)/(k+1)$ for every $Q\in \cal F$, for every $Q\in \cal F$ there is a positive integer $j\le k$ such that $x_j$ belongs to $\pi_\la(Q)$. Let $j(Q)$ be the smallest such $j$. For $1\le j\le k$, let $m_j$ be the number of rectangles $Q$ of $\cal F$ such that $j(Q)=j$. Then we have $\sum_{j=1}^{k} m_j^2 \le N_{\delta}(\la)$, and so, if $B=\{1\le j\le k \mid m_j\ne 0\}$, $\sum_{j\in B} m_j^2\le N_{\delta}(\la)$. We have also $\sum_{j\in B}m_j=\sum_{j=1}^k m_j=\#{\cal F}$. 

Since 
$$\# B\cdot N_{\delta}(\la)\ge \# B\sum_{j\in B}m_j^2\ge(\sum_{j\in B}m_j)^2=(\#{\cal F})^2,$$
we have 
$$\# B\ge \frac{(\#{\cal F})^2}{N_{\delta}(\la)}\ge \frac{b^2\cdot \delta^{-2(d+d')}}{\tilde c \cdot \delta^{-d-d'}}={\tilde c}^{-1}\cdot b^2 \cdot \delta^{-(d+d')}.$$
Given $j \in B$ we may choose $Q \in {\cal F}$ such that $j(Q)=j$, and an interval $J(j)$ such that $x_j \in J(j)\subset \pi_\la(Q)$ with size $(R+1)/(k+1)$. The interiors of these intervals are all disjoint
and so  
$$Leb(\pi_{\la}(X))\ge \frac{\# B\cdot (R+1)}{(k+1)}\ge \frac{\# B\cdot (R+1)}{2k}\ge \frac{\#B\cdot c_0^{-1}\delta}4\ge \frac{{\tilde c}^{-1} b^2 \delta^{-(d+d')}c_0^{-1}\delta}4=\frac{c_0^{-1} {\tilde c}^{-1} b^2\delta^{1-d-d'}}4.$$

\vskip 0,5cm
{\bf 4 - Dimension estimates - constructing a recurrent set of good scales}.
\vglue .2in

{\bf 4.1} 
We define sets $E(\underline
\theta , \underline \theta')$ as the sets of $s \in J_R$ which verify:
$$N_{\rho^{1/m}}(\underline \theta , \underline \theta', s) \leq c_5 \rho^{-(d+d')/m}
\; ,
\leqno(1)$$
$$\displaystyle \sum_{(\underline b, \underline b')} N_{\rho^{1/m}}(T_{\underline
b}T'_{\underline b'} (\underline \theta , \underline \theta', s))
\leq c_5 \; \hat
\rho^{-(d+d')}\rho^{-(d+d')/m}\leqno(2)$$ where, in (2), $\hat \rho$ takes the values
$\rho^{1/m} ,
\rho^{2/m} , \dots , \rho^{(m-1)/m}$ and
$(\underline b, \underline b')$ describes the pairs of $\Sigma(\hat \rho) \times
\Sigma'(\hat
\rho)$ beginning by $(\theta_0 , \theta'_0)$. We then put
$$E(\underline a, \underline a')= \displaystyle \bigcup_{(\underline \theta,
\underline
\theta')}  E(\underline \theta, \underline \theta') \; ,\leqno (3)$$

\noindent where $\; \underline a \in \Sigma(\rho), \underline a' \in
\Sigma'(\rho)$ and
$(\underline
\theta, \underline \theta')$ describes the pairs of $\Sigma^- \times \Sigma^{'-}$
ending by
$(\underline a, \underline a')$. The integer $m$ is fixed once for all, and
should satisfy $m \geq 3 \; , \; m > \alpha^{-1}$. The real number $\alpha \in (0,1]$ denotes here the
parameter intervening in the scale recurrence lemma.

We remark that the hypothesis 2.2, 2.4 which allow us to prove the scale recurrence lemma are less restrictive when the parameter $\alpha$ is smaller. This means that we may always choose the parameter $\alpha$ satisfying $\alpha \le 1/2$, which we will do in what follows. 

We assume that the scale recurrence property holds (which, as we have seen, is generically true).

Now, if we choose $c_5$ large enough, Proposition 1 from section 3 guarantees that the condition
$$ \hbox {Leb} (J_R \setminus E(\underline a, \underline a')) \leq c_1 \leqno (4)$$

\noindent
which appears in the scale recurrence lemma holds. So the scale recurrence lemma gives us, for $\underline a \in \Sigma(\rho),
\underline a'
\in \Sigma'(\rho)$ a family of compact sets $E^*(\underline a, \underline a')
\subset J_R$ which satisfy the conclusions (i), (ii) and (iii) of 2.3.

We will consider the compact set $\tilde {\cal L} \subset {\cal S}$ constituted by the triples
$(\underline \theta ,
\underline \theta', s)$ such that $ s \in E^*(\underline
a , \underline a')$, where $\underline \theta , \underline
\theta'$ finish respectively by $\underline a \in \Sigma(\rho) \; , \; \underline a'
\in \Sigma'(\rho)$  .

\vskip 0,5 cm
{\bf 4.2} 

The fiber $L(\underline \theta , \underline \theta', s) \subset {\re}$ (for
$(\underline
\theta ,
\underline \theta',s) \in \tilde {\cal L})$ is constituted by the values of $t$ for which there is at least one pair  $(\underline b^1 , \underline b^{'1}) \in
\Sigma(\rho) \times \Sigma'(\rho)$ which (putting $T_{\underline b^1} \; T'_{\underline b^{'1}}
(\underline
\theta ,
\underline \theta', s)=(\underline \theta^1 , \underline \theta^{'1} , s_1))$ has the following properties:
$$\hbox{for } \vert \tilde s_1 - s_1 \vert \leq \frac{3}{4} \; \rho^\alpha,
\hbox{we have }
(\underline
\theta^1 , \underline \theta^{'1},  \tilde s_1) \in \tilde {\cal L}  \;  \; ; \leqno
(7)$$
$$t \hbox{ belongs to the image
of the rectangle }
I(\underline b^1)
\times I'(\underline b^{'1}) \hbox{ by } \pi_{\underline \theta , \underline \theta',s} \; .\leqno (8)$$

In this chapter we will prove the following

{\bf Proposition} - There is $c_6 > 0$ such that we have
$$Leb (L (\underline \theta , \underline \theta' , s)) \geq c_6 \rho^{1-(d+d')}$$ for every
$(\underline
\theta , \underline \theta' , s) \in \tilde {\cal L} \; .$

\vskip 0,5 cm
{\bf 4.3 - Proof of the proposition}

{\bf 4.3.1} Let us fix $(\underline
\theta , \underline \theta' , s) \in \tilde {\cal L}$. By definition of
$\tilde {\cal L}$, there
exists $(\underline a , \underline a') \in \Sigma(\rho) \times \Sigma'(\rho)$
ending
$(\underline
\theta , \underline \theta')$ such that $s \in E^*(\underline a , \underline a')$.
By the
conclusion (i) of the scale recurrence lemma, there exists $\tilde s
\in E(\underline a ,
\underline a')$ such that
$$\vert s - \tilde s \vert \leq c_2 \; \rho^\alpha .\leqno(9)$$

By the definition (3) of $E(\underline a , \underline a')$, there is
$(\underline {\tilde
\theta} ,
\underline {\tilde \theta}\,') \in  \Sigma^- \times \Sigma^{'-}$ ending by
$(\underline a
,
\underline a')$ such that
$$\tilde s \in E(\underline{{\tilde \theta}} , {\underline{\tilde \theta}}\,').
 \leqno(10)$$

The above relations (1), (2) are therefore verifyied by $(\underline{{\tilde
\theta}}
,\underline{{\tilde \theta}} \,' , \tilde s)$.

From the conclusion (iii) of the scale recurrence lemma, there are at least
$c_3 \;
\rho^{-(d+d')}$ pairs $(\underline b , \underline b') \in \Sigma(\rho) \times
\Sigma'(\rho)$
(with
$b_0 = \theta_0  \; , \; b'_0 = \theta'_0)$ such that, if we put
$$T_{\underline b} T'_{\underline b'} (\underline  \theta ,\underline
\theta\,', s)=
(\underline {\hat \theta} , \underline {\hat \theta}\,', \hat s) \; , $$
we have $[\hat s -
\rho^\alpha , \hat s + \rho^\alpha] \subset E^*(\underline b , \underline b')$.
We thus have $(\underline {\hat \theta} , \underline {\hat
\theta}\,',
\hat s_1) \in \tilde {\cal L}$ for $\vert \hat s - \hat s_1 \vert \leq
\rho^\alpha$.

We will denote by $\Gamma_m$ the set of these pairs $(\underline b , \underline b')$.
We thus have
$$\# \Gamma_m \geq c_3 \; \rho^{-(d+d')} \; .\leqno(11)$$

{\bf 4.3.2 -} Since the map $\underline \theta \to k^{\underline \theta}$ is not necessarily Lipschitz, but is $\alpha$-Hölder (with $\alpha>1/m$), and the sets $E^*(\underline a , \underline a')$ are contained in neighbourhoods of $E(\underline a , \underline a')$ of sizes $O(\rho^\alpha)$ (but not $O(\rho)$), we need to consider a sequence of scales in a geometric progression of ratio $\rho^{1/m}$ before aplying the Scale Recurrence Lemma: We now define, by descending recursion on $1
\leq \ell \leq m$, a subset $\Gamma_\ell$ of $\Sigma(\rho^{\ell/m}) \times
\Sigma'(\rho^{\ell/m})\; .$ We just defined $\Gamma_m$. For $\ell < m$, the definition depends on a 
large enough constant $c'_5$:
$\Gamma_\ell$ is constituted of the pairs $(\underline b , \underline b')$
satisfying
$$N_{\rho^{1/m}}(T_{\underline
b}T'_{\underline b'} (\underline{{\tilde \theta}},\underline{{\tilde \theta}} \,' , \tilde s)) \leq c'_5 
\rho^{-(d+d')/m} \;
;\leqno(12)$$
$$\hbox{there are at least} \quad c^{'-1}_5 \rho^{-\frac{1}{m}(d+d')} \quad
\hbox{pairs in }
 \Gamma_{\ell+1} \hbox{ beginning by }  (\underline b , \underline
b')\; . \leqno(13)$$

From the relation (2) (for $(\underline {\tilde \theta} , \underline
{\tilde
\theta} ' ,
\tilde s))$, there are at most $cc_5 \, c^{'-1}_5 \; \rho^{-\ell (d+d')/m}$
pairs in
$\Sigma(\rho^{\ell/m}) \times \Sigma'(\rho^{\ell/m})$ which do not satisfy
(12). On the other hand, there are at most $c \rho^{-(d+d')/m}$ pairs in $\Sigma(\rho^{(\ell
+1)/m}) \times
\Sigma'(\rho^{(\ell+1)/m})$ extending a given pair of
$\Sigma(\rho^{\ell/m}) \times
\Sigma'(\rho^{\ell/m})$. We conclude that if $c'_5$ is large enough, we will have by
descending recursion from (11) :
$$\# \Gamma_\ell \geq c_7 \; \rho^{-\frac{\ell}{m} (d+d')} \; .\leqno(14)$$

{\bf 4.3.3 -} For words $\underline c , \underline c'$ of $\Sigma ,
\Sigma'$ beginning
respectively by $\theta_0 , \theta'_0$, let us denote by $J(\underline c , \underline
c')$ (resp.$\tilde J(\underline c , \underline c'))$ the interval of ${\re}$ image
of the rectangle
$I(\underline c)
\times I'(\underline c')$ by $\pi_{\underline \theta , \underline \theta',s}$ (resp. by
$\pi_{\underline {\tilde \theta} , \underline {\tilde \theta}',\tilde s}$).

For $1 \leq \ell <\ell' \leq m \; , (\underline b , \underline b') \in
\Gamma_\ell$, let us denote by
$\Gamma_{\ell'} (\underline b , \underline b')$ the set of pairs
$(\underline c ,
\underline c') \in \Gamma_{\ell'}$ which begin by $(\underline b ,
\underline b')$ ;
let us denote by $K(\underline b , \underline b') \subset J(\underline b , \underline b')$
the union over
$\Gamma_m(\underline b ,
\underline b')$ of the intervals $J(\underline c , \underline c')$.

{\bf Lemma :} There is $c_8 > 0$ such that we have, for $1 \leq \ell \leq m \;
, \;
(\underline b ,
\underline b') \in \Gamma_\ell$ :

$$Leb(K (\underline b , \underline b')) \geq c_8\rho^{(m-\ell) (1-d-d')/m} Leb(J(\underline b ,
\underline b'))$$

{\bf 4.3.4 -} The two following elementary results are used in the proof of the lemma and proved in
[MY, p. 70].

\noindent {\bf Sublemma 1 -} Let $(J_\alpha)_{\alpha \in A} ,
(J'_\alpha)_{\alpha \in A}$
two families of intervals and $\varepsilon>0 , \lambda	>1 , \nu >0$ constants satisfying:

\noindent (i) $\qquad$ For every $\alpha \in A$, we have $\varepsilon < Leb(J_\alpha) <
\lambda
\; \varepsilon$ and $\varepsilon < Leb(J'_\alpha)$.

\noindent (ii) $\qquad$ For every $\alpha \in A$, the centers of $J_\alpha ,
J'_\alpha$ are
distant from at most $\nu \varepsilon$ ;

We then have
$$Leb \; (\bigcup J'_\alpha) \geq (\lambda(4 \nu +4))^{-1} \; Leb \; (\bigcup
J_\alpha)$$

\noindent {\bf Sublemma 2 -} Let $(J_\alpha)_{\alpha \in A}$ be a family
of intervals,
and $0 < \nu <1$ ; for each $\alpha \in A$, let $K_\alpha$ a subset of
$J_\alpha$
satisfying $Leb(K_\alpha) \geq
\nu Leb(J_\alpha)$. We have then
$$Leb (\bigcup_\alpha \; K_\alpha) \geq \frac{1}{2} \; \nu \; Leb
(\bigcup_\alpha \;
J_\alpha)\; . $$

{\bf 4.3.5 - Proof of the Lemma}

It will be done by descending recursion on $\ell$ ; the case $\ell = m$ is
tautologic (with
$c_8=c_8(m)=1)$ and so we assume $\ell <m$.

By Proposition 2 from section 3, for each subset $X \subset K \times
K'$ which is the intersection of $K \times K'$ with a disjoint union of at least $b\cdot \rho^{(-d-d')/m}$ rectangles of the type $I_{\und a} \times I_{\und a'}'$, with $(\und a, \und a') \in \Sigma(\rho^{1/m}) \times \Sigma'(\rho^{1/m})$, if $N_{\rho^{1/m}}(\underline {\tilde \theta} , \underline {\tilde \theta}', \tilde s) \le c_5 \rho^{-(d+d')/m}$, we have 
$$Leb \; (\pi_{\underline {\tilde \theta} , \underline {\tilde \theta}', \tilde
s} (X)) \geq
c_0^{-1} \cdot {\tilde c}^{-1} \cdot b^2\cdot\rho^{(1-d-d')/m}/4. $$ 

Let $(\underline b \, , \, \underline b')\in \Gamma_\ell$. We apply
the preceding inequality replacing $({\underline {\tilde \theta} , \underline {\tilde
\theta}',
\tilde s})$ by $T_{\underline b} T_{\underline b'}\;({\underline {\tilde
\theta} ,
\underline{\tilde \theta}',
\tilde s})$ and taking for $X$ the version renormalized by $T_{\underline
b} \;
T'_{\underline b'}$ of the union of the rectangles $I(\underline c) \times
I(\underline c')$,
where $(\underline c \; , \underline c')$ describes $\Gamma_{\ell +1}(\underline
b \; , \;
\underline b')$. From (12) and (13), and using the preceding inequality  (with $b=c^{'-1}_5$), we get

$$Leb(\bigcup_{\Gamma_{\ell +1}(\underline b \; , \; \underline b')}
\tilde J(
\underline c ,
\underline c')) \geq c_9 \rho^{(1-d-d')/m} Leb(\tilde J (\underline b , \underline b')) \;
.\leqno(15)$$

We will now compare the measures of  $\cup \tilde J(\underline c , \underline
c')$ and $\cup
J(\underline c , \underline c')$ using sublemma 1. Since we have
$d(\underline \theta
, \tilde {\underline \theta)} \leq c \rho$, the $C^1$-distance between
$k^{\underline \theta}$ and
$k^{\tilde {\underline \theta}}$ is at most of the order of $\rho^\alpha$.

The same holds for $k'\,^{\underline \theta'}$ and $k'\,^{\underline{\tilde
\theta}'}\, $. We
also have $\vert s - \tilde s \vert \leq c_2 \, \rho^\alpha$. The $C^1$-distance
between
$\pi_{\underline \theta, \underline \theta',s}$ and
$\pi_{\underline {\tilde \theta}, \underline {\tilde \theta}', \tilde s}$ is
thus at most of the order of
$\rho^{\alpha}$. We recall that $\alpha > \frac{1}{m}$. Sublemma 1 (with
$\varepsilon
\sim \rho^{1/m}$) and (15) give therefore:
$$Leb(\bigcup_{\Gamma_{\ell +1}(\underline b \; , \; \underline b')} J(
\underline c ,
\underline c')) \geq c \; Leb(\bigcup_{\Gamma_{\ell +1}(\underline b \; , \;
\underline b')} \tilde J(
\underline c ,
\underline c')) \geq c c_9 \rho^{(1-d-d')/m}  Leb(\tilde J (\underline b , \underline b')) 
\geq  \leqno(16)$$
$$\geq c c'_9 \rho^{(1-d-d')/m}  Leb(J (\underline b , \underline b')).$$
By the recursion hypothesis, for $(\underline c , \underline c') \in
\Gamma_{\ell
+1}(\underline b \; ,
\underline b') \; , $ we have
$$Leb(K(\underline c , \underline c')) \geq c_8(\ell+1) \rho^{(m-\ell-1)(1-d-d')/m} Leb(J
(\underline c ,
\underline c')) \; , \leqno(17)$$ which gives us, after (16) and sublemma 2 :
$$Leb(\bigcup_{\Gamma_{\ell +1}(\underline b \; , \; \underline b')} K(
\underline c ,
\underline c')) \geq \frac{1}{2} \; c_8(\ell +1) \rho^{(m-\ell-1)(1-d-d')/m} \; Leb(\bigcup_{\Gamma_{\ell
+1}(\underline b
\; , \; \underline b')} J( \underline c ,
\underline c')) \geq  \leqno(18)$$
$$\frac{1}{2} \; c'_9 \rho^{(1-d-d')/m}  \; c_8 (\ell+1) \rho^{(m-\ell-1)(1-d-d')/m} \; Leb(J
(\underline b ,
\underline b')) =  \frac{1}{2} \; c'_9 \; c_8 (\ell+1) \rho^{(m-\ell)(1-d-d')/m} \; Leb(J
(\underline b ,
\underline b')) \;  ,$$
which is the inequality of the lemma with $c_8(\ell)= \frac{1}{ 2} \; c'_9
\; c_8(\ell + 1) \; . $ \qed

\vglue .2in
{\bf 4.3.6 - End of the proof of the Proposition}

As in the proof of the Lemma, we apply Proposition 2 from section 3, now to the union of the rectangles of $\Gamma_1$. Since $\# \Gamma_1 \geq c_7 \rho^{-\frac{1}{m} (d+d')}$, the union for $(b, b')\in \Gamma_1$ of $\tilde J(b, b')$ has Lebesgue measure at least $c_0^{-1} \cdot {\tilde c}^{-1} \cdot {c_7}^2\cdot\rho^{(1-d-d')/m}/4$. Since, by the Lemma, for every $(b, b')\in \Gamma_1$, $Leb(K (\underline b , \underline b')) \geq c_8\rho^{(m-1) (1-d-d')/m} Leb(J(\underline b ,\underline b'))$, using sublemmas 1 and 2 as in the proof of the previous Lemma we conclude that the union for $(b, b')\in \Gamma_1$ of $K (\underline b , \underline b')$ has Lebesgue measure at least $c_6 \rho^{1-(d+d')}$, where $c_6>0$ is a constant, and this implies the Proposition. \qed

\vskip 0,5cm
{\bf 5 - The main results}.

\noindent{\bf Proposition}: Let $K$\,, $K'$ be regular Cantor sets with $HD(K)+HD(K')<1$
such that $K$ satisfies property $H_{\alpha}$ for some $\alpha>0$ as in the preceeding discussion. Suppose that, for any $(\underline \theta , \underline \theta' , s) \in \Sigma^- \times \Sigma^{'-} \times J_R$ there are admissible finite words $\underline
b, \underline b'$ such that $T_{\underline
b}T'_{\underline b'} (\underline \theta , \underline \theta', s) \in \tilde {\cal L}$. Then, given a $C^1$ map $h$ from a neighbourhood of $K \times K'$ to $\re$ such that, in some point of $K \times K'$ its gradient is not parallel to any of the two coordinate axis, we have $HD(h(K \times K'))=d+d'$.
\vglue .1in

\noindent{\bf Proof}: 
The inequality $HD(h(K \times K')) \le d+d'$ follows from the fact that $h$ is Lipschitz and $HD(K \times K')=d+d'$.

By hypothesis, and by continuity of $dh$, we may find a pair of periodic points $p$ and $p'$ of $K$ and $K'$, respectively, with addresses  $\und{\ov a}=\und a \und a \und a...$ and  $\und{\ov a'}=\und a' \und a' \und a'...$, where $\und a$ and $\und a'$ are finite sequences, such that $dh(p,p')$ is not a real multiple of $dx$ nor of $dy$. There are increasing sequences of natural numbers $(m_k),(n_k)$ and $s \in J_R$ such that $$\lim_{k \to \infty}-\frac{\frac{\partial h}{\partial y}(p,p').|I'_{{\und a'}^{n_k}}|}{\frac{\partial h}{\partial x}(p,p').|I_{{\und a}^{m_k}}|}=s.$$ 
Since $h$ is $C^1$, for large $k$, the restriction of $h$ to $I_{{\und a'}^{n_k}} \times I'_{{\und a}^{m_k}}$ is very close in the $C^1$ topology, up to a composition on the left by an affine map, to the linear map $x-sy$. So, the composition of $h$ with the map $(f^{{\und a'}^{n_k}},f^{'{\und a'}^{n_k}})$, which brings $I(a_0) \times I'(a'_0)$ to $I_{{\und a}^{m_k}} \times I'_{{\und a'}^{n_k}}$ (here $f^{{\und a'}^{n_k}},f^{'{\und a'}^{n_k}}$ are inverse branches of iterates of the expansive dynamics which define $K$ and $K'$) is very close in the $C^1$ topology, up to a composition on the left with an affine map, to $\pi_{\und{\ov a} , \und{\ov a'}, s}$.

By hypothesis, there is $\delta>0$ and finite admissible words $\underline b, \underline b'$ such that $T_{\underline b}T'_{\underline b'} (\und{\ov a} , \und{\ov a'}, s) \in \tilde {\cal L}$, so the composition of $h$ with the map $(k^{{\und a'}^{n_k} \und b},k^{'{\und a'}^{n_k} \und b'})$, which brings $I(a_0) \times I'(a'_0)$ to $I_{{\und a}^{m_k} \und b} \times I'_{{\und a'}^{n_k} \und b'}$ is very close in the $C^1$ topology, up to a composition on the left with an affine map, to $\pi_{T_{\underline b}T'_{\underline b'} (\und{\ov a} , \und{\ov a'}, s)}$.  

Given $\eta>0$, suppose that $\rho$ is so small that $\frac{c_7}{2 c_0} (\frac{c^{-1}_0}{2+2R})^{d+d'-\eta}>\rho^{\eta}$.
For $(\underline \theta , \underline \theta', s') \in \tilde {\cal L}$, the Proposition above implies that there are intervals which are image
of rectangles of the type $I(\underline b^1) \times I'(\underline b^{'1})$ (with $(\underline b^1 , \underline b^{'1}) \in
\Sigma(\rho) \times \Sigma'(\rho)$ ) by  $\pi_{\underline \theta , \underline \theta',s}$ such that, putting $T_{\underline b^1} \; T'_{\underline b^{'1}}(\underline \theta , \underline \theta', s)=(\underline \theta^1 , \underline \theta^{'1} , s_1)$, for  $\vert \tilde s_1 - s_1 \vert \leq \frac{3}{4} \; \rho^\alpha$, we have $(\underline \theta^1 , \underline \theta^{'1},  \tilde s_1) \in \tilde {\cal L}$, and such that the measure of the union of these intervals is at least $c_7 \rho^{1-(d+d')}$. We may select a disjoint union of these intervals whose measure is at least $\frac 12 c_7 \rho^{1-(d+d')}$. Since for words $\underline a$ of $\Sigma(\rho)$ we have $c^{-1}_0 \rho \leq \vert I( \underline a ) \vert \leq c_0 \rho $, we have at least  $\frac 12 c_0^{-1}c_7 \rho^{-(d+d')}$ intervals in the above disjoint union, which are contained in the interval $\pi_{\underline \theta , \underline \theta',s}(I(\theta_0)\times I'(\theta_0'))$, which has size at most $2R+2$. So, the sum of the $(d+d'-\eta)$-powers of the ratios between the sizes of these intervals and the size of the interval $\pi_{\underline \theta , \underline \theta',s}(I(\theta_0)\times I'(\theta_0'))$ (which contains their union) is at least 
$\frac 12 c_0^{-1}c_7 \rho^{-(d+d')} (\frac{c_0^{-1} \rho}{2+2R})^{d+d'-\eta}=\frac{c_7}{2 c_0} (\frac{c^{-1}_0}{2+2R})^{d+d'-\eta}  \rho^{-\eta}>1$, by hypothesis. These facts constitute an open phenomenon in the $C^1$ topology, in the sense that, if we change $\pi_{\underline \theta , \underline \theta',s}$ by a very close $C^1$ map $F$, the conclusions remain true. Moreover, this phenomenon is hereditary: in this context, if $F$ is very close to $\pi_{\underline \theta , \underline \theta',s}$, then the restriction of $F$ to a rectangle $I(\underline b^1) \times I'(\underline b^{'1})$ as above is even closer in the $C^1$ topology, up to composing on the left with an affine map, to a map $\pi_{\underline \theta^1 , \underline \theta^{'1},  \tilde s_1}$, with $(\underline \theta^1 , \underline \theta^{'1},  \tilde s_1) \in \tilde {\cal L}$ (since limit geometries behave as attractors for renormalization operators, and differentiable maps are very close to linear maps in small scales).   
 
Now, the conclusion of the proposition is a consequence of the following fact, mentioned in [MY2], which may be proved with usual techniques of fractal geometry (as the mass distribution principle; see section 4.1 of [F]):

Let us consider, for $r \geq 1$, a finite family $\mathcal{J}(r)$ of disjoint compact intervals. We suppose that each interval of $\mathcal {J}(r + 1)$ is contained in an interval of $\mathcal{J}(r)$ and that we have

$$ \sup_{ r \geq 1} \; \sup_{\substack {I' \subset  I \\ I \subset \; \mathcal{J}(r),\; I' \in \; \mathcal {J}(r + 1)}} \frac{|I|}{|I'|} \ < +\infty$$

If we have, for every $r \geq 1$, $I \in \mathcal{J}(r)$:

$$ \sum_{\substack{I'\subset \; I \\ I' \in \; \mathcal{J}(r + 1)}} {\left(\frac{|I'|}{|I|}\right )}^{\tilde d} \geq 1\;,$$

then we have 

$$HD(\;\bigcap_{r \geq 1}\bigcup_{\mathcal{J}(r)}I\;) \geq \tilde {d} \; \text .$$

Indeed, the previous discussion implies that we may apply this statement to the set \linebreak $h(K({\und a}^{m_k}) \times K'({\und a'}^{n_k}))$, for $k$ large, with $\tilde d=d+d'-\eta$. Since $\eta>0$ is arbitrary, the result follows.
\qed

\vglue .2in
\noindent{\bf Corollary 1}: Let $K$\,, $K'$ be regular Cantor sets such that $HD(K)+HD(K')<1$
satisfying properties $H_{\alpha}$ for some $\alpha>0$ then:

i) If $K$ and $K'$ have periodic orbits whose ratio of logarithms of the norms of their eigenvalues is irrational, then, given a $C^1$ map $h$ from a neighbourhood of $K \times K'$ to $\re$ such that, in some point of $K \times K'$ its gradient is not parallel to any of the two coordinate axis, we have $HD(h(K \times K'))=HD(K) + HD(K')$. Consequently, if $K$ has two periodic orbits whose ratio of logarithms of the norms of their eigenvalues is irrational, then, given a $C^1$ map $h$ from a neighbourhood of $K \times K'$ to $\re$ such that, in some point of $K \times K'$ its gradient is not parallel to any of the two coordinate axis, we have $HD(h(K \times K'))=HD(K) + HD(K')$.

ii) If $K$ has a sequence of periodic points $(x_n)$ with corresponding eigenvalues $(\mu_n)$ and $K'$ has a periodic point $y$ with eigenvalue $\lambda$ such that $(\log |\mu_n|/\log |\lambda|)$ is a sequence of rational numbers whose denominators tend to $+\infty$, then, given a $C^1$ map $h$ from a neighbourhood of $K \times K'$ to $\re$ such that, in some point of $K \times K'$ its gradient is not parallel to any of the two coordinate axis, we have $HD(h(K \times K'))=HD(K) + HD(K')$. Consequently, if $K$ has a sequence of periodic points $(x_n)$ with corresponding eigenvalues $(\mu_n)$ such that $(\log |\mu_n|/\log |\mu_1|)$ is a sequence of rational numbers whose denominators tend to $+\infty$, then, given a $C^1$ map $h$ from a neighbourhood of $K \times K'$ to $\re$ such that, in some point of $K \times K'$ its gradient is not parallel to any of the two coordinate axis, we have $HD(h(K \times K'))=HD(K) + HD(K')$. 

\noindent{\bf Proof}: We will prove that, under conditions i) or under conditions ii), then for any $(\underline \theta , \underline \theta' , s) \in \Sigma^- \times \Sigma^{'-} \times J_R$ there are admissible finite words $\underline
b, \underline b'$ such that $T_{\underline
b}T'_{\underline b'} (\underline \theta , \underline \theta', s) \in \tilde {\cal L}$, and so the result follows from the previous Proposition.

Let $N$ be the set of $(\underline \theta , \underline \theta' , s) \in \Sigma^- \times \Sigma^{'-} \times J_R$ for which there are admissible finite words $\underline
b, \underline b'$ such that $T_{\underline
b}T'_{\underline b'} (\underline \theta , \underline \theta', s) \in \tilde {\cal L}$. Since $\tilde{\cal L}$ is nonempty (and, by the conclusion (ii) of the scale recurrence lemma 2.3, meets both orientations $s>0$ and $s<0$), we deduce that $N$ is nonempty, and even that for any given $(\und{\te}, \und{\te'}) \in \Sigma^- \times \Sigma^{\prime -}$ there are $\la_1 < 0 < \la_2$ such that $(\und{\te}, \und{\te'}, \la_i) \in N$ for $i=1,2$: indeed, given $(\underline a , \underline a')
\in \Sigma(\rho) \times \Sigma'(\rho)$ such that $E^*(\underline a , \underline a')$ is nonempty and meets both orientations, given $(\und{\te}, \und{\te'}) \in \Sigma^- \times \Sigma^{\prime -}$ we take finite words $\underline d, \underline d'$ larger than $\underline a , \underline a'$ and ending by $\underline a , \underline a'$, respectively, and, since $T_{\und{d}}T'_{\und{d}'}(\und{\te},\und{\te}',s)=(\und{\te}\sigma(\und{d}),\und{\te}'\sigma(\und{d}'),\lambda s)$, for some $\lambda \neq 0$, and $\sigma(\und{d}), \sigma(\und{d})$ end respectively by $\underline a , \underline a'$, we may chose $\la_1 < 0 < \la_2$ for which the conclusion holds. By compacity of $\Sigma^- \times \Sigma^{\prime -}$, there is $\delta>0$ such that, given $(\underline b, \underline b')\in \Sigma(\delta) \times \Sigma'(\delta)$, there are $\la_1(\underline b, \underline b') < 0 < \la_2(\underline b, \underline b')$ such that if $\und{\te}\in\Sigma^-$ and $\und{\te'} \in \Sigma^{\prime -}$ end respectively by $\underline b$ and $\underline b'$, then, for every $\la \in (\la_1(\underline b, \underline b')-\delta,\la_1(\underline b, \underline b')+\delta)\cup(\la_2(\underline b, \underline b')-\delta,\la_2(\underline b, \underline b')+\delta)$, we have $(\und{\te}, \und{\te'}, \la) \in N$.

Suppose that there are  periodic points $x$, $x'$ of period $n$, $n'$ of $g$, $g'$ respectively such that
$$
\frac{\log |Dg^n(x)|}{\log |Dg^{\prime n'}(x')|}
$$
is irrational. Let $\und{\ov c}=\und c \und c \und c...$ and $\und{\ov c'}=\und c' \und c' \und c'...$be the addresses of $x$ and $x'$, respectively, where $\und c$ and $\und c'$ are finite sequences of sizes $n$ and $n'$ ending by $\tilde c$ and $\tilde c'$, respectively. 
Given $(\und{\te},\und{\te}',s)\in \Sigma^- \times \Sigma^{\prime -} \times \re^*$, we have, for $m, r$ large positive integers, 
$$T_{\tilde c\und{c}^m}T'_{\tilde c'{\und{c}'}^r}(\und{\te},\und{\te}',s)=(\und{\te}\und{c}^m,\und{\te}'{\und{c}'}^r,\ve \ve' s|I^{\und{\te}}(\tilde c\und{c}^m)|^{-1}|I^{\und{\te}'}(\tilde c'{\und{c}'}^r)|),$$ where where $\ve = +1$ (resp. $-1$) if $f_{\tilde c\und{c}^m}$ is orientation preserving (resp.
reversing) and $\ve' = +1$ (resp. $-1$) if $f'_{\tilde c'{\und{c}'}^n}$ is orientation preserving (resp.
reversing).
Since $(\und{\te}\und{c}^m,\und{\te}'{\und{c}'}^r)$ converges to $(\und{\ov c},\und{\ov c'})$ and 
$|I^{\und{\te}}(\tilde c\und{c}^m)|^{-1}|I^{\und{\te}'}(\tilde c'{\und{c}'}^r)|$ is asymptotically $\tilde \lambda |Dg^{\prime n'}(x')|^r/|Dg^n(x)|^m$, for some constant $\tilde \lambda>0$, the hypothesis on the eigenvalues imply that we may choose $m, r$ in such a way that, if $\und{\ov c},\und{\ov c'}$ end respectively by $\underline b \in \Sigma(\delta)$ and $\underline b'\in \Sigma'(\delta)$, 
then $$\ve \ve' s|I^{\und{\te}}(\tilde c\und{c}^m)|^{-1}|I^{\und{\te}'}(\tilde c'{\und{c}'}^r)|) \in 
(\la_1(\underline b, \underline b')-\delta,\la_1(\underline b, \underline b')+\delta)\cup(\la_2(\underline b, \underline b')-\delta,\la_2(\underline b, \underline b')+\delta).$$
This implies that $T_{\tilde c\und{c}^m}T'_{\tilde c'{\und{c}'}^r}(\und{\te},\und{\te}',s)\in N$, and, since $(\und{\te},\und{\te}',s)\in \Sigma^- \times \Sigma^{\prime -} \times \re^*$ is arbitrary, we finally get that $N = \Sigma^- \times \Sigma^{\prime -} \times \re^*$.  This proves i).

In order to prove ii), if $\log |\mu_n|/\log |\lambda|=p/q$, with $gcd(p,q)=1$, then there are positive integers $r, s$ such that $pr-qs=1$, and so $log(|\mu_n|^r/|\lambda|^s)=r\log |\mu_n|/-s\log |\lambda|=\log |\lambda|/q$. So, for any positive integer $k$, $\log(|\mu_n|^{rk}/|\lambda|^{sk})=k\log |\lambda|/q$. If $q$ is large enough (depending on $\delta$), these ratios become sufficiently dense in $(0,+\infty)$ in such a way that the preceding argument implies that, for any $(\underline \theta , \underline \theta' , s) \in \Sigma^- \times \Sigma^{'-} \times J_R$, there are admissible finite words $\underline
b, \underline b'$ such that $T_{\underline
b}T'_{\underline b'} (\underline \theta , \underline \theta', s) \in \tilde {\cal L}$. 

If $K$ has a sequence of periodic points $(x_n)$ with corresponding eigenvalues $(\mu_n)$ such that $(\log |\mu_n|/\log |\mu_1|)$ is a sequence of rational numbers whose denominators tend to $+\infty$, choose $y$ a periodic point of $K'$ and $\lambda$  its eigenvalue. If $(\log |\mu_1|/\log |\lambda|)$ is irrational, we conclude by i) that, given a $C^1$ map $h$ from a neighbourhood of $K \times K'$ to $\re$ such that, in some point of $K \times K'$ its gradient is not parallel to any of the two coordinate axis, we have $HD(h(K \times K'))=HD(K) + HD(K')$. If $(\log |\mu_1|/\log |\lambda|)$ is rational, then the ratios $(\log |\mu_n|/\log |\lambda|)=(\log |\mu_n|/\log |\mu_1|)\cdot (\log |\mu_1|/\log |\lambda|)$ form a sequence of rational numbers whose denominators tend to $+\infty$, and  we conclude by ii) that, given a $C^1$ map $h$ from a neighbourhood of $K \times K'$ to $\re$ such that, in some point of $K \times K'$ its gradient is not parallel to any of the two coordinate axis, we have $HD(h(K \times K'))=HD(K) + HD(K')$.
\qed

\vglue .15in
In order to remove the hypothesis $HD(K)+HD(K')<1$, we need the following Lemma.

\noindent{\bf Lemma :} Given a regular Cantor set $K$ defined by an expanding map $g$ and constants $a, b$ with $0\le a<b\le HD(K)$, there is a regular Cantor set $\tilde K$ contained in $K$ satisfying $a<HD(\tilde K)<b$.

\noindent{\bf Proof}: Let ${\Sigma}^{(n)}$ be the set of words of $\Sigma$ of length $n$ and ${\mathcal{R}}^{(n)}=\{I_{\underline{b}}, \underline{b}\in {\Sigma}^{(n)}\}$. Fix $\tilde{c}, \tilde{d} \in {\bf A}$ such that $(\tilde{d},\tilde{c})\in {\cal B}$. Let $X^n=\{\underline{b}=(b_1,\cdots,b_n)\in {\Sigma}^{(n)}: b_1=\tilde{c}, b_n=\tilde{d}\}$. For any positive integer $r$, and $\underline{b_1}, \underline{b_2},\dots,\underline{b_r} \in X^n$, we have $\underline{b_1} \ \underline{b_2}\dots\underline{b_r} \in X^{nr}\subset {\Sigma}^{(nr)}$. Let ${\mathcal{\tilde{R}}}^{(n)}=\{I_{\underline{b}}, \underline{b}\in X^n\}$.

For $R\in {\mathcal{R}}^{(n)}$ take $\Lambda_{n,R}=\displaystyle\sup \left|(\psi^{n})^{\prime}_{|_{R}}\right|$ and $\lambda_{n,R}=\displaystyle\inf \left|(\psi^{n})^{\prime}_{|_{R}}\right|$. By the bounded distortion property and the mixing condition, there is ${\hat c}>0$ such that ${\hat c}\Lambda_{n,R} \le \lambda_{n,R} \le \Lambda_{n,R}, \forall R\in {\mathcal{R}}^{(n)}$ and
$${\hat c}\sum_{R\in\mathcal{R}^{(n)}}(\Lambda_{n,R})^{-d}\le \sum_{R\in\mathcal{\tilde{R}}^{(n)}}(\Lambda_{n,R})^{-d}\le \sum_{R\in\mathcal{R}^{(n)}}(\Lambda_{n,R})^{-d}, \forall d\in [0,2], n\ge 1.$$
On the other hand, from [PT], pp. 69-70, it follows that, if we define $d_n$ implicitely by
$$\sum_{R\in\mathcal{R}^{(n)}}(\Lambda_{n,R})^{-d_n}=1,$$
then $\lim d_n=HD(K)$, so in particular for $n$ large we have $d_n>(a+b)/2$. Notice also that there is $\lambda_1>1$ such that $\Lambda_{n,R}\ge \lambda_1^n, \forall n\ge 1$.

Let $\varepsilon=(b-a)/2>0$ and let $n$ be large so that $d_n>HD(K)-\varepsilon/2$ and $\lambda_1^{n\varepsilon}>2/{\hat c}$. We have
$$\sum_{R\in\mathcal{\tilde{R}}^{(n)}}(\Lambda_{n,R})^{-(a+b)/2}\ge {\hat c}\sum_{R\in\mathcal{R}^{(n)}}(\Lambda_{n,R})^{-(a+b)/2}>{\hat c}.$$
Removing elements of $\mathcal{\tilde{R}}^{(n)}$ one by one, we get a subset ${\mathcal D}^{(n)}$ of $\mathcal{\tilde{R}}^{(n)}$ (which is equal to $\{I_{\underline{b}}, \underline{b}\in Y^n\}$, for a certain subset $Y^n$ of $X^n$) and a word $\underline{\tilde c}\in X^n\setminus Y^n$ such that 
$$\sum_{R\in \mathcal{D}^{(n)}\cup \{I_{\underline{\tilde c}}\}}(\Lambda_{n,R})^{-(a+b)/2}>{\hat c} \text{ and } \sum_{R\in\mathcal{D}^{(n)}}(\Lambda_{n,R})^{-(a+b)/2}\le {\hat c}.$$
Notice that 
$$\sum_{R\in\mathcal{D}^{(n)}}(\Lambda_{n,R})^{-(a+b)/2}\ge \sum_{R\in \mathcal{D}^{(n)}\cup \{I_{\underline{\tilde c}}\}}(\Lambda_{n,R})^{-(a+b)/2}-\lambda_1^{-n (a+b)/2}>{\hat c}/2.$$
and so $$\sum_{R\in\mathcal{D}^{(n)}}(\Lambda_{n,R})^{-a}=\sum_{R\in\mathcal{D}^{(n)}}(\Lambda_{n,R})^{\varepsilon-(a+b)/2}>1.$$
Moreover, we have
$$\sum_{R\in\mathcal{D}^{(n)}}(\lambda_{n,R})^{-(a+b)/2}\le \sum_{R\in\mathcal{D}^{(n)}}({\hat c}\Lambda_{n,R})^{-(a+b)/2}<{\hat c}^{-1}{\hat c}=1.$$

\noindent We may define the regular Cantor set (with expanding map $g^n$) $$\displaystyle\tilde{K}:=\bigcap_{r\ge 0} g^{-nr}(\cup_{\hat I \in \mathcal{D}^{(n)}}\hat I).$$  
The previous estimates imply that $$\sum_{R\in{\mathcal{D}}^{(n)r}}(\Lambda_{nr,R})^{-a}\ge \left(\sum_{R\in\mathcal{D}^{(n)}}(\Lambda_{n,R})^{-a}\right)^r>1,$$
where $\mathcal{D}^{(n)r}=\{I_{\underline{{\hat c}} \underline{c_2}\dots\underline{c_r}}, \underline{c_j}\in Y^n, \forall j\le r\}$.
On the other hand, we have
$$\sum_{R\in{\mathcal{D}}^{(n)r}}(\lambda_{nr,R})^{-(a+b)/2}\le \left(\sum_{R\in\mathcal{D}^{(n)}}(\lambda_{n,R})^{-(a+b)/2}\right)^r<1,$$
an it follows from the estimates of [PT], pp. 69-70 that $a<HD(\tilde K)<(a+b)/2<b$. 
\qed

\vglue .15in
\noindent{\bf Corollary 2}: Let $K$\,, $K'$ be regular Cantor sets
such that $K$ satisfies property $H_{\alpha}$ for some $\alpha>0$. Suppose that

i) $K$ has two periodic orbits whose ratio of logarithms of the norms of their eigenvalues is irrational

or

ii) $K$ has a sequence of periodic points $(x_n)$ with corresponding eigenvalues $(\mu_n)$ such that $(\log |\mu_n|/\log |\mu_1|)$ is a sequence of rational numbers whose denominators tend to $+\infty$.

Then, given a $C^1$ map $h$ from a neighbourhood of $K \times K'$ to $\re$ such that, in some point of $K \times K'$ its gradient is not parallel to any of the two coordinate axis, we have $HD(h(K \times K'))=\min\{1,d+d'\}$.
\vglue .1in

\noindent{\bf Proof}: If $d+d'<1$ it follows directly from Corollary 1. If $d+d'\ge 1$, given $\delta>0$, it follows from the Lemma above that there is a regular Cantor set $K''$ contained in $K'$ satisfying $1-d-\delta<HD(K'')<1-d\le d'$. Using Corollary 1 we conclude that $1-\delta<HD(h(K \times K'))\le 1$ for every $\delta>0$, so, in this case, $HD(h(K \times K'))=1$, which implies the result.
\qed
\vskip.1in
\noindent{\bf Remark} - In the case when $K$ and $K'$ are regular Cantor sets such that $HD(K)+HD(K')>1$, it was proved in Theorem 2 of [MR] using the results of [MY] that, under open and dense conditions in the $C^r$-topology, for any $r>1$ (namely if $(g, g')$ belong to the set $V$ defined in [MY], where $g$ and $g'$ are the expanding maps which define $K$ and $K'$), for any $C^1$ map $f$ from a neighbourhood of $K \times K'$ to $\re$ such that, in some point of $K \times K'$ its gradient is not parallel to any of the two coordinate axis, $f(K \times K')$ has non-empty interior in $\re$.
\vglue .2in
\noindent{\bf Corollary 3}: Let $K$ and $K'$ be, respectively, a stable and an unstable Cantor set of a horseshoe $\Lambda$. Suppose that $K$ satisfies property $H_{\alpha}$ for some $\alpha>0$ and

i) $K$ has two periodic orbits whose ratio of logarithms of the norms of their eigenvalues is irrational

or

ii) $K$ has a sequence of periodic points $(x_n)$ with corresponding eigenvalues $(\mu_n)$ such that $(\log |\mu_n|/\log |\mu_1|)$ is a sequence of rational numbers whose denominators tend to $+\infty$.

Then, given a $C^1$ map $h$ from a neighbourhood of $\Lambda$ to $\re$ such that, in some point of $\Lambda$ its gradient is not parallel to the stable nor to the unstable direction of $\Lambda$ at this point, then we have $HD(h(\Lambda))=\min\{1,HD(\Lambda)\}$.

\noindent{\bf Proof}: It follows from Corollary 2 and from the fact that, given a point of a horseshoe $\Lambda$ in a surface, there is a $C^1$ diffeomorphism from a neighbourhood $\cal U$ of it to an open rectangle of ${\mathbb R}^2$ with sides parallel to the coordinate axis such that the image of local stable leaves of $\Lambda$ are contained in horizontal lines, the image of local unstable leaves of $\Lambda$ are contained in vertical lines and the image of $\Lambda \cap {\cal U}$ is a cartesian product ${\tilde K} \times \tilde{K'}$ of differentiable embeddings ${\tilde K}$ and $\tilde{K'}$ of $K$ and $K'$, respectively (notice that $HD(\Lambda)=HD(K)+HD(K')$).     
\qed

\vglue .15in

\noindent{\bf Proof of Theorem 1:}  We want to check the hypothesis of Corollary 2. Let $h\colon \Sigma^+ \to K$ be the homeomorphism which conjugates the dynamics of $g$ over $K$ to the unilateral sub-shift $\sigma:\Sigma^+\to \Sigma^+$. We have by hypothesis two limit geometries $\und \alpha$ and $\und \beta$ in $\Sigma^-$ which end by the same letter $a_0$, and a point which will be identified by its address $\und \gamma$ in $\Sigma$, which begins by $a_0$, such that $(k_{\und \alpha} \circ k_{\und \beta}^{-1})''(k_{\und \beta}(h(\und \gamma)) \neq 0$. If we change $\und \alpha$ and $\und \beta$ by other very close limit geometries (whose last letters are very similar to those of $\und \alpha$ and $\und \beta$), and $h(\und \gamma)$ by another very close point in $K(a_0)$ (whose address begins in the same way as $\gamma$), we will have the corresponding second derivative still non-zero, and of the same sign. So, there exist, by continuity, finite words $\und a, \und b$, which we may assume that begin by the same $c_0$ such that $a_0 c_0$ is an allowed transition, which end in a very similar way as $\und \alpha$ and $\und \beta$, respectively, and a finite word $\und c=a_0 \und c'$, which ends by $a_0$ and begins in a very similar way as $\und \gamma$, but not in a very similar way as $\und a$, such that, for every infinite sequences $\und x, \und y, \und z$ we have $(k_{\und y \und a}\circ k_{\und z \und b}^{-1})''(k_{\und z \und b}(h(\und c \und x))$ non-zero.

    Consider now the limit geometry associated to the sequence $\und{\ov a}=\und a \und a \und a...$, and the $|\und a|-$th iterated $F$ of the dynamics which defines $K$ restricted to the interval $I(a_0 \und a)$. Then, by the remarks at the end of section 1.2, the map $\tilde F:=k_{\und{\ov a}} \circ F \circ (k_{\und{\ov a}})^{-1}$, conjugated to $F$, is affine, and brings the interval $I^{\ov a}(a_0 \und a)$ on $I(a_0)$. Its derivative $\lambda$ is the eigenvalue of the periodic orbit of $K$ of period $|\und a|$ which address is $\und{\ov a}$. We will call by $G$ the $|\und b|-$th iterate of the dynamics which defines $K$ restricted to the interval $I(a_0 \und b)$, and by $\tilde G$ the map $k_{\und{\ov a}} \circ G \circ (k_{\und{\ov a}})^{-1}$. By the remarks at the end of section 1.2, $\tilde G$ is equal to $k_{\und{\ov a}} \circ k_{\und{\ov a} \und b}^{-1}\circ L$, where $L$ is an affine map taking $I^{\ov a}(a_0 \und b)$ to $I(a_0)$, and so $(\tilde G)''$ restricted to $I^{\ov a}(a_0 \und b \und c')$ does not vanish. Let $H$ be the  $|\und c'|-$th iterated of the dynamics which defines $K$ restricted to the interval $I(\und c)$, and $\tilde H:=k_{\und{\ov a}} \circ H \circ (k_{\und{\ov a}})^{-1}$. If we put $\tilde R=\tilde H \circ (\tilde G|_{I(a_0 \und b \und c')})$, he have that $\tilde R$ is equal to $k_{\und{\ov a}} \circ k_{\und{\ov a} \und b \und c'}^{-1}\circ {\tilde L}$, where ${\tilde L}$ is an affine map taking $I^{\ov a}(a_0 \und b \und c')$ to $I(a_0)$. 
		
Let us now consider the two following pre-periodic points of $K$: $u=h(\und c \und{\ov a})$ and $v=h(a_0 \und b \und c' \und{\ov a})$; let us put $\tilde u=k_{\und{\ov a}}(u)$ and $\tilde v=k_{\und{\ov a}}(v)$. We have $G(v)=u$, and so $\tilde G (\tilde v)=\tilde u$. Since $(\tilde G)''(\tilde v) \neq 0$, at least one of the following possibilities holds, by the chain rule: $(\tilde H)''(\tilde u) \neq 0$ or $(\tilde R)''(\tilde v) \neq 0$. Indeed, we have $(\tilde R)'=(\tilde H)' \circ \tilde G \cdot (\tilde G)'$ and $(\tilde R)''=(\tilde H)'' \circ \tilde G \cdot (\tilde G)'^2+(\tilde H)' \circ \tilde G \cdot (\tilde G)''$, and so $(\tilde R)''(\tilde v)=(\tilde H)''(\tilde u) \cdot ((\tilde G)'(\tilde v))^2+(\tilde H)'(\tilde u) \cdot (\tilde G)''(\tilde v)$. In the first case, the sequence of the periodic points $p_k$ whose addresses have periods $\und c {\und a}^k$, $k\ge 1$ is such that the sequence of logarithms of the eigenvalues of their orbits modulo $\log \lambda$ is convergent and has eventually distinct terms. Indeed, conjugating by $k_{\und{\ov a}}$, we get the sequence of points $\tilde p_k=k_{\und{\ov a}}(p_k)$ such that $(\tilde G)^k\circ \tilde H(\tilde p_k)=\tilde p_k$, and the eigenvalue of the orbit of $p_k$ coincides with $((\tilde G)^k\circ \tilde H)'(\tilde p_k)=\lambda^k (\tilde H)'(\tilde p_k)$. Since the points $(\tilde p_k)$ are distinct and converge to $\tilde u=k_{\und{\ov a}}(u)$, our claim follows from $(\tilde H)''(\tilde u) \neq 0$. In the second case, the same thing happens for the sequence of periodic points $q_k$ whose addresses have periods $a_0 \und b \und c' {\und a}^k$, $k\ge 1$ (and the proof is analogous to the previous case). In any case, we obtain sequences of pairs of periodic points (one of them, in each pair, being $h(\und{\ov a})$), whose ratio of logarithms of the norms of the eigenvalues are either irrational or rational with arbitrarily large denominators, which solves the problem, by Corollary 1. \qed

\vskip 0,5cm

\end{document}